\documentclass[12pt,a4paper,oneside]{amsart}
\usepackage[a4paper]{geometry}
\usepackage{mathpazo} 

\usepackage{mathrsfs}
\usepackage[UglyObsolete]{diagrams} 

\usepackage[vcentermath,noautoscale]{youngtab}
\Yboxdim{10pt}  

\usepackage[colorlinks=false]{hyperref}

\DeclareMathOperator{\ad}{\mathsf{ad}}
\DeclareMathOperator{\Der}{Der}
\DeclareMathOperator{\Homol}{H}
\DeclareMathOperator{\ess}{\mathsf{S}}

\newcommand{\vertbar}{\>|\>}
\newcommand{\set}[2]{\ensuremath{\{ #1 \vertbar #2 \}}}

\newtheorem{question}{Question}

\begin{document}

\title{Lie algebras and around: selected questions}
\author{Pasha Zusmanovich}
\address{Department of Mathematics \\
University of Ostrava \\
30. dubna 22, 70103 Ostrava \\
Czech \mbox{Republic}}
\email{pasha.zusmanovich@osu.cz}
\date{First written April 25, 2016; last revised May 7, 2016}
\keywords{
Lie algebra; Zassenhaus algebra; cohomology; deformations; Koszul dual operads;
sum of subalgebras; characteristic $2$; Ado Theorem; Whitehead Lemma;
ternary maps
}
\subjclass[2010]{17B05; 17B10; 17B50; 17B55; 17B56; 05A99}
\dedicatory{
To my teacher and friend Askar Serkulovich Dzhumadil'daev on his 60th birthday}

\begin{abstract}
Several open questions are discussed. The topics include cohomology of current 
and related Lie algebras, algebras represented as the sum of subalgebras,
structures and phenomena peculiar to characteristic $2$, and variations on 
themes of Ado, Whitehead, and Banach.
\end{abstract}

\maketitle

\section*{Introduction}

I am presenting here a, perhaps, haphazard collection of questions I am 
interested in. Being haphazard as it is, this collection features somewhat 
unexpected and, hopefully, fascinating connections between different topics.

This is a modest contribution in honor of Askar Dzhumadil'daev. I first met him
in 1983, when, as an undergraduate student, I started to participate in his
seminar on Lie algebras. Since then and throughout many years, I enjoyed his 
wisdom, unfailing enthusiasm, friendship, and support. Most of what I know in 
mathematics I owe to him.

\section{Cohomology of Lie algebras coming from Koszul dual operads}\label{sec-young}

Current Lie algebras -- that is, Lie algebras of the form $L \otimes A$, where
$L$ is a Lie algebra and $A$ is a commutative associative algebra -- as well as
algebras close to them, are ubiquitous in mathematics and physics (sufficient is to mention Lie algebras of smooth functions on a manifold 
prominent in gauge theory, Kac--Moody Lie algebras, modular semisimple Lie algebras, etc.). The Lie bracket is defined in an obvious 
factor-wise way:
$$
[x \otimes a, y \otimes b] = [x,y] \otimes ab ,
$$
where $x,y\in L$, $a,b\in A$. A vast generalization of this construction comes 
from the operad theory: if $A$ is an algebra over a binary quadratic operad 
$\mathcal P$, and $B$ is an algebra over the operad Koszul dual to $\mathcal P$,
then the tensor product $A \otimes B$ equipped with the bracket
\begin{equation}\label{eq-brack}
[a \otimes b, a^\prime \otimes b^\prime] = 
aa^\prime \otimes bb^\prime - a^\prime a \otimes b^\prime b ,
\end{equation}
where $a,a^\prime \in A$, $b,b^\prime \in B$, is a Lie algebra.

Due to a big flexibility of this construction, many interesting Lie algebras 
can be represented in this form, for a suitable pair of Koszul dual operads
and algebras over them. Perhaps the most remarkable examples, beyond current
Lie algebras, are various algebras appearing in physics 
(Schr\"odinger--Virasoro, Heisenberg--Virasoro, etc.): they are represented
as tensor products of algebras over left and right Novikov operads. This 
remarkable observation was implicit already in the pioneering works of 
I. Gelfand and S.P. Novikov and their collaborators (\cite{gelfand-dorfman} and
\cite{novikov}); after that, Pei and Bai (\cite{pei-bai} and references therein)
noted that Lie algebras in question admit realization as affinizations of certain left 
Novikov algebras; Dzhumadil'daev noted in \cite{dzhu-novikov} that left and right Novikov operads are Koszul dual to each other; the explicit claim
was made in \cite[\S 5]{vavilov-fest} by putting all the pieces together.

Therefore, it seems to be interesting to develop a general method for computing
cohomology and other invariants of Lie algebras given by bracket 
(\ref{eq-brack}). By this, we mean to express cohomology or other invariants
of such Lie algebras in terms of invariants of tensor factors $A$ and $B$, 
similarly how it was done earlier for low-degree cohomology of current Lie
algebras.

Let us try to pursue further an approach for computation of (co)homology of 
current Lie algebras described in \cite[\S 4]{low}. (For simplicity, we will 
consider cohomology with coefficients in the trivial module $K$, assume some 
finiteness conditions, and zero characteristic of the base field, but similar 
considerations may be applied in more general settings). Let us decompose the
modules $\bigwedge^n (A \otimes B)$ according to the well-known
Cauchy formula:
$$
\bigwedge\nolimits^n (A \otimes B) \simeq \bigoplus_{\lambda \vdash n} 
\Big(Y_\lambda(A) \otimes Y_{\lambda^\sim}(B)\Big) .
$$
Here the direct sum of vector spaces is taken over all Young diagrams of size 
$n$, $\lambda^\sim$ is the Young diagram transposed to $\lambda$, and 
$Y_\lambda$ is the Young symmetrizer corresponding to $\lambda$. Assuming that
at least one of the algebras $A, B$ is finite-dimensional, passing to the decomposition
of the dual vector spaces, and decomposing the differential 
$d: (\bigwedge^n(A \otimes B))^* \to (\bigwedge^{n+1}(A \otimes B))^*$ in the 
Chevalley--Eilenberg complex accordingly, we get the following maps on the
Young graph:

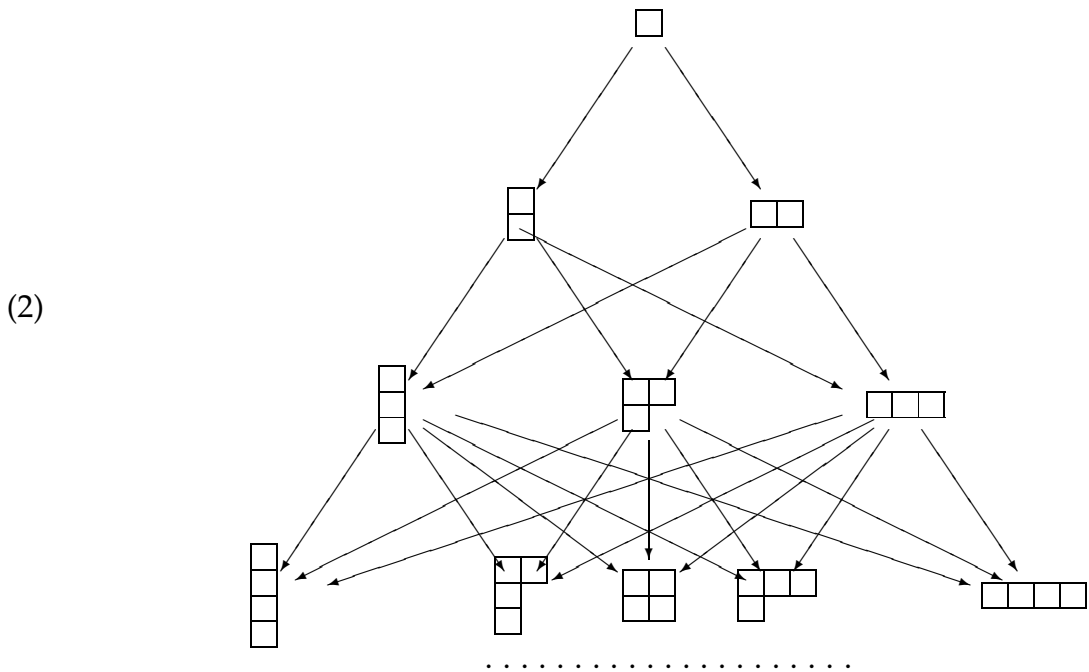
\begin{equation}\label{eq-young}
\begin{diagram}[width=2em,height=3em]
&&&&&& \yng(1) &&&&&&&                                                 \\ 
&&&&& \ldTo && \rdTo &&&&&                                             \\ 
&&&& \yng(1,1) &&&& \yng(2) &&&&                                       \\ 
&&& \ldTo && \rdTo\rdTo(5,2) & & \ldTo(6,2)\ldTo && \rdTo &&&          \\ 
&& \yng(1,1,1) &&&& \yng(2,1) &&&& \yng(3) &&                          \\ 
& \ldTo && \rdTo\rdTo(4,2)\rdTo(6,2)\rdTo(10,2) && \ldTo(6,2)\ldTo & \dTo & \rdTo\rdTo(6,2) && \ldTo(10,2)\ldTo(6,2)\ldTo(4,2)\ldTo && \rdTo &  \\ 
\yng(1,1,1,1) &&&& \yng(2,1,1) && \yng(2,2) && \yng(3,1) &&&& \yng(4)  \\ 
\end{diagram}
\end{equation}
\begin{center}
. . . . . . . . . . . . . . . . . . . . .
\end{center}
\bigskip
Here each Young diagram $\lambda$ denotes the vector space 
$Y_\lambda(A)^* \otimes Y_{\lambda^\sim}(B)^*$, and arrows denote the 
corresponding components of the differential $d$.

Intuitively it should be clear that the more arrows here vanish, the easier it
would be to compute the corresponding cohomology. In the case of the pair of 
operads (Lie, associative commutative), i.e. for current Lie algebras 
$L \otimes A$, a miracle happens: approximately half of the arrows vanish (roughly, those going from ``left'' to 
``right''), what allows to define a certain spectral sequence on the 
Chevalley--Eilenberg complex computing the cohomology 
$\Homol^*(L \otimes A, K)$. In the low cohomology degrees and/or for particular
types of algebras, this spectral sequence allows to express 
$\Homol^*(L \otimes A, K)$ in terms of cohomology and other invariants of the
tensor factors $L$ and $A$. Unfortunately, this miracle fails for the other
pairs of Koszul dual operads, even such a classical one as 
(associative, associative).

\begin{question}
What makes the pair (Lie, associative commutative) special in this regard?
In which other situations (i.e., for a particular pair of Koszul dual operads,
or for a particular kind of algebras over some pair of Koszul dual operads)
``many'' arrows in the corresponding Young graph (\ref{eq-young}) will vanish? 
In particular, for which types of Lie algebras expressed as the tensor products
of left Novikov and right Novikov algebras, this will happen?
\end{question}

\section{Algebras represented as the sum of subalgebras}

My mathematical debut, under the guidance of Askar Dzhumadil'daev, was the proof
of the Kegel--Kostrikin conjecture about solvability of a modular 
finite-dimensional Lie algebra $L$ represented as the vector space sum 
$L = N + M$ of two nilpotent subalgebras $N$, $M$ (\cite{mat-zametki}; around 
the same time this was established also by Panyukov, \cite{pan1}). The statement
is true over fields of characteristic $p>2$, but in characteristic $2$ there is
a counterexample found by Petravchuk, \cite{pet}. Take the 
$3$-di\-men\-si\-o\-nal characteristic $2$ analog of $\mathsf{sl}_2$: the simple Lie algebra $\ess$, with the basis $\{e_{-1}, e_0, e_1\}$ subject to multiplication 
$$
[e_{-1},e_0] = e_{-1}, \quad [e_{-1},e_1] = e_0, \quad [e_0,e_1] = e_1 .
$$
Its $2$-envelope $\ess^{[2]}$ is $5$-dimensional and admits the decomposition 
$\ess^{[2]} = N \oplus M$, where the $2$-dimensional abelian subalgebra $N$ is 
linearly spanned by $e_0 + e_{-1} + e_{-1}^{[2]}$ and $e_0 + e_1 + e_1^{[2]}$, 
and the $3$-dimensional nilpotent subalgebra $M$ is linearly spanned by 
$e_{-1}^{[2]}$, $e_0$ and $e_1^{[2]}$ (the vector space sum in this case is 
direct).

Do such examples exist in higher dimensions? Of course, any 
(fi\-ni\-te- or in\-fi\-ni\-te-di\-men\-si\-o\-nal) current Lie algebra 
$\ess \otimes A$ admits such a decomposition:
$$
\ess \otimes A = (N \otimes A) \oplus (M \otimes A) ,
$$ 
so it provides such an example provided it itself is non-solvable 
(for example, when $A$ contains a unit).

A slightly more interesting example can be obtained as an extension of the
corresponding current Lie algebras of the form $\ess \otimes A + \mathcal D$,
where $\mathcal D$ acts on $A$ by derivations, what includes semisimple Lie 
algebras. Namely, we have the decomposition
\begin{equation}\label{eq-decomp}
\ess \otimes A + \mathcal D = 
(N \otimes A + \mathcal D) \oplus (M \otimes A) .
\end{equation}

An easy induction on $n$ proves that for a Lie algebra 
$L = S \otimes A + \mathcal D$ with $A$ unital, and for any positive integer 
$n$, we have
$$
L^n = S^n \otimes A + 
      \sum_{\substack{i+j=n \\ i>1, j\ge 1}} S^i \otimes A \mathcal D^j(A) +
      S \otimes \mathcal D^{n-1}(A) + \mathcal D^n .
$$
This implies that if $N$ is nilpotent, and $\mathcal D$ is nilpotent as an 
algebra of derivations of $A$ (and hence is nilpotent as an abstract Lie 
algebra), then the algebra $N \otimes A + \mathcal D$ is nilpotent too. 
Therefore, (\ref{eq-decomp}) provides a decomposition of a nonsolvable Lie 
algebra into the sum of nilpotent subalgebras.

Yet it would be more interesting to generalize Petravchuk's decomposition for 
an arbitrary Zassenhaus algebra $W_1^\prime(n)$ in characteristic $2$. 
Zassenhaus algebras appear prominently in ongoing efforts of classification of 
simple Lie algebras in characteristic $2$ (cf. 
\cite[Vol. I, \S 7.6]{strade-intro}, \cite{S}, \cite{grishkov-zus}, and 
references therein). In characteristic $p=2$, unlike for $p>2$, the Zassenhaus algebra has dimension $2^n-1$ and can be defined as the algebra with the basis $\set{e_i}{-1 \le i \le 2^n-3}$
subject to multiplication
$$
[e_i,e_j] = \begin{cases}
\binom{i+j+2}{i+1} e_{i+j} &\text{if } -1 \le i+j \le 2^n-2 \\
0                          &\text{otherwise}.
\end{cases}
$$

The algebra $\ess = W_1^\prime(2)$ is the first algebra in the series. The
$2$-envelope $W_1^\prime(n)^{[2]}$ of $W_1^\prime(n)$ coincides with the 
derivation algebra of $W_1^\prime(n)$, has dimension $2^n + n - 1$, and is 
spanned, in addition to elements of $W_1^\prime(n)$, by elements 
$(\ad e_{-1})^{2^k}$, $k = 1, 2, \dots, n-1$, and $(\ad e_{2^{n-1}-1})^2$.

\begin{question}
Find a link with combinatorial interpretation of the number $2^n + n - 1$ 
as the shortest length of a sequence of $0$ and $1$ containing all subsequences 
of length $n$ (see \cite[A052944]{eis}).
\end{question}

\begin{question}
Is it true that $W_1^\prime(n)^{[2]}$ admits a decomposition into the sum of 
two nilpotent subalgebras?
\end{question}

Virtually nothing is known about the Kegel--Kostrikin question in the 
infinite-dimensional case -- beyond almost obvious cases when one imposes some
sort of finiteness conditions on one or both of the summands; all such cases are
reduced quickly to the finite-dimensional situation.

As a first step, one may wish to prove that such an algebra satisfies a 
nontrivial (Lie) identity. According to \cite[Corollary 2.5]{ultra}, a Lie 
algebra $L$ does not satisfy a nontrivial identity if and only if a free Lie 
algebra is embedded into an ultraproduct of $L$. As the ultraproduct 
construction obviously commutes with the decomposition into the sum of 
subalgebras, the question whether a Lie algebra $L = N + M$ does not satisfy a 
nontrivial identity is equivalent to the question whether its ultraproduct 
$L^{\mathscr U} = N^{\mathscr U} + M^{\mathscr U}$ does not contain a free Lie
subalgebra. As being nilpotent (of a fixed nilpotency index) is the first-order
property, by the {\L}o\'s theorem the Lie algebras $N^{\mathscr U}$ and 
$M^{\mathscr U}$ are also nilpotent. Thus the question whether the sum of two
nilpotent Lie algebras satisfies a nontrivial identity, is equivalent to
an apriori more special 

\begin{question}
Is it true that an infinite-dimensional Lie algebra represented as the sum
of two nilpotent subalgebras, cannot contain a free Lie algebra as a subalgebra?
\end{question}

In the theory of (associative) PI algebras, there is a similar long-standing 
open question: whether the sum of two PI algebras is PI? (See \cite{klm} and 
(numerous) references therein). By the same reasonings as in the Lie case, this
question is equivalent to

\begin{question}
Is it true that an associative algebra represented as the sum of two PI 
subalgebras, cannot contain a free associative algebra as a subalgebra?
\end{question}

Another interesting topic is to develop a machinery to express the (co)homology
(Chevalley--Eilenberg, Hochschild, etc.) of such algebras in terms of factors 
and their action on each other. In the case when the sum of subalgebras is 
direct, one may attempt to mimic the approach of \S \ref{sec-young}, albeit
in a more simple situation, as we have direct sums instead of tensor products.
If say, a Lie algebra $L = M \oplus N$ is represented as the vector space direct
sum of subalgebras $M$ and $N$, then, instead of the Cauchy formula we have a 
more simple isomorphism of vector spaces:
$$
\bigwedge\nolimits^n (N \oplus M) \simeq 
\bigoplus_{\substack{i+j=n \\ i,j \ge 0}} 
\bigwedge\nolimits^i(N) \otimes \bigwedge\nolimits^j(M) .
$$

Passing to the dual vector spaces, and decomposing the differentials in the
Chevalley--Eilenberg complex, as in \S \ref{sec-young}, we get the following
picture (now instead of the Young graph we have a more simpler triangle):
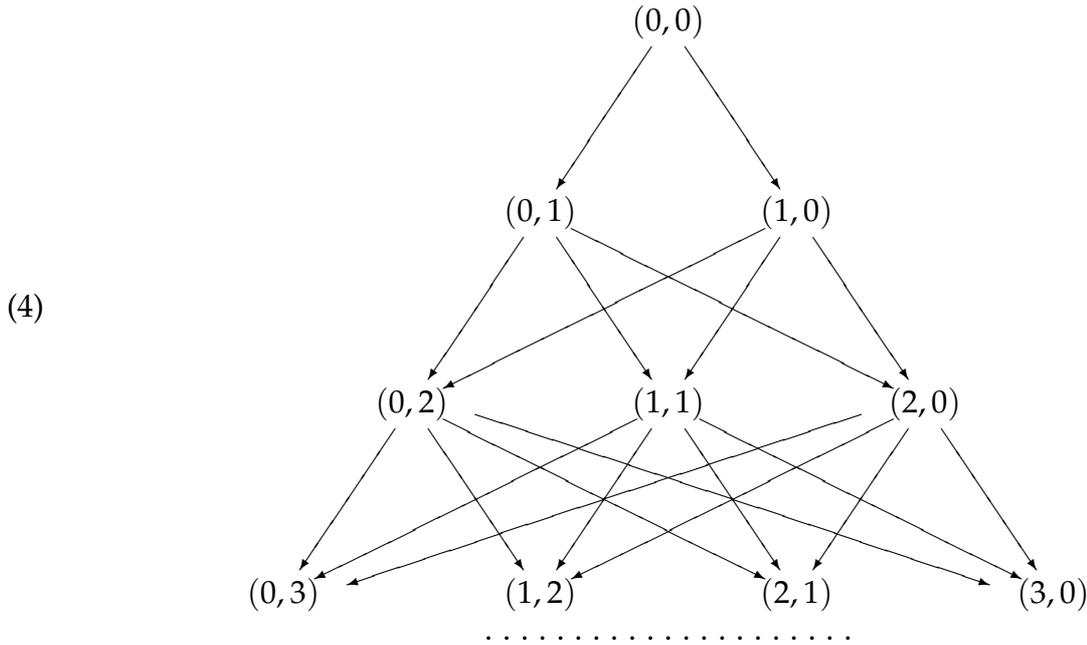
\begin{equation}\label{eq-wedge}
\begin{diagram}[width=2em,height=3em]
&&&&&& (0,0) &&&&&                                             \\ 
&&&&& \ldTo && \rdTo &&&                                       \\ 
&&&& (0,1) &&&& (1,0) &&                                       \\ 
&&& \ldTo && \rdTo\rdTo(6,2) && \ldTo(6,2)\ldTo && \rdTo &     \\ 
&& (0,2) &&&& (1,1) &&&& (2,0) &                               \\ 
& \ldTo && \rdTo\rdTo(6,2)\rdTo(10,2) && \ldTo(6,2)\ldTo && \rdTo\rdTo(6,2) && \ldTo(10,2)\ldTo(6,2)\ldTo && \rdTo &                           \\ 
(0,3) &&&& (1,2) &&&& (2,1) &&&& (3,0)                         \\ 
\end{diagram}
\end{equation}
\begin{center}
. . . . . . . . . . . . . . . . . . . . .
\end{center}
\bigskip
Here each pair $(i,j)$ denotes the vector space 
$\bigwedge^i(N)^* \otimes \bigwedge^j(M)^*$.

\begin{question}
Are there any patterns (vanishing or otherwise) in the graph (\ref{eq-wedge})
in the general situation? In some special cases?
\end{question}

Ideally, a positive answer to this question should allow to develop a 
cohomological machinery which would unify and generalize various situations: 
some particular instances of the Hochschild--Serre spectral sequence, 
a stuff related to ``Tate Lie algebras'', ``Japanese cocycles'' (see, e.g., 
\cite[\S 2.7]{chiral}), etc.

A somewhat similar machinery is contained in an interesting and seemingly 
entirely forgotten paper \cite{dupre} (there, the author presents an alternative
derivation of the Lyndon--Hochschild--Serre spectral sequence for the semidirect
product of groups, but similar considerations seem to be applicable as well to 
the group-theoretic analog of our situation, i.e. for a group $G = AB$ 
decomposed into the product of its subgroups $A$ and $B$; the promised sequels 
to \cite{dupre} treating the Lie-algebraic and associative cases have never
appeared).

\section{
Deformations and ``commutative'' cohomology in characteristic 
\texorpdfstring{$2$}{2}
}

Classification of finite-dimensional simple Lie algebras over algebraically 
closed fields of characteristic zero is a classical piece of mathematics, 
crystallized at the end of XIX--beginning of XX centuries. It took the mankind
another some 100 years to achieve the same classification over fields of 
characteristic $p>3$ (see \cite{strade-intro}). The cases $p=2$ and $3$ remain 
widely open. In \cite{grishkov-zus}, an attempt was made to advance the case 
$p=2$ basing on earlier results of Skryabin \cite{S}. The general line of attack
is more or less the same as in ``big'' characteristics: one first classifies 
algebras of small toral rank, and then, taking advantage of appropriate root 
space decompositions, glue the results together; also, many questions are 
reduced to computation of deformations of certain classes of algebras. In 
particular, in \cite{S} simple Lie algebras having a Cartan subalgebra of toral
rank $1$ are characterized as certain filtered deformations of semisimple Lie algebras $L$ such that
\begin{equation}\label{eq-semi}
S \otimes \mathcal O_1(n) \subset L \subseteq 
\Der(S) \otimes \mathcal O_1(n) + K\partial ,
\end{equation}
where $\mathcal O_1(n)$ is the divided powers algebra, $\partial$ its standard
derivation, and either $n=2$ and $S \simeq W_1^\prime(n)$, or $n=1$ and $S$ is
isomorphic to a two-variable Hamiltonian algebra. 

In \cite{grishkov-zus}, these deformations were computed in the simplest case
$S \simeq \ess$, what allowed to classify simple Lie algebras of absolute toral
rank $2$ and having a Cartan subalgebras of toral rank $1$ -- a small but 
necessary step in the classification program. To further advance along this 
road, one need to compute these deformations in all the cases.

\begin{question}\label{q-deform}
Compute deformations of semisimple Lie algebras in characteristic $2$ of the
form (\ref{eq-semi}).
\end{question}

In the process of these computations, it became apparent that a new cohomology
theory peculiar to characteristic $2$ plays a role. This cohomology is defined
via the standard formula for the coboundary map in the Chevalley--Eilenberg
complex, with the alternating cochains being replaced by symmetric ones. Note
that we can (profitably) consider commutative $2$-cocycles of Lie algebras in 
arbitrary characteristic (\cite{alia-d} and \cite{comm2}), albeit they do not 
lead to any cohomology; while in characteristic $2$ we have a bona fide 
cohomology theory. Unlike the usual cohomology, this complex is apriori not 
restricted by the dimension of the algebra, so new interesting phenomena, 
similar to those appearing in cohomology of Lie superalgebras (in any 
characteristic), occur. Generally, this ``commutative'' cohomology is different 
from the Chevalley--Eilenberg cohomology. For example, while the second 
cohomology of the Zassenhaus algebra $W_1^\prime(n)$ with coefficients in the 
trivial module is zero (note that this is in striking difference with the cases
of ``big'' characteristics; if $p>3$, the corresponding cohomology is 
$1$-dimensional, leading to the modular analog of the famous Virasoro algebra, 
cf. \cite{h2-zass}), the analogous ``commutative'' cohomology has dimension $n$
and is generated by ``commutative'' $2$-cocycles
\begin{equation*}
e_i \vee e_j \mapsto \begin{cases}
1 &\text{if } i=j=2^k - 2, \text{or } \{i,j\} = \{-1,2^{k+1}-3\} \\
0 &\text{otherwise}.
\end{cases}
\end{equation*}
for $k = 0,\dots,{n-1}$.

(This can be established by considering subalgebras of $W_1^\prime(n)$ spanned
by a ``small'' number of root vectors -- what corresponds to the cases $n=2$ and
$3$ -- similarly how it was done in computation of commutative $2$-cocycles
on simple Lie algebras of classical type in \cite{alia-d}).

Besides a few isolated computations like just presented, virtually nothing is 
known about this kind of cohomology, so any result about it would be of 
interest. For example, to compute deformations in Question \ref{q-deform}, one need to compute 
low-degree ``commutative'' cohomology with various coefficients of simple Lie 
algebras involved -- the Zassenhaus and Hamiltonian algebras.

\begin{question}
Compute the ``commutative'' cohomology for various Lie algebras in 
characteristic $2$.
\end{question}

\begin{question}
Is it possible to represent the ``commutative'' cohomology as a derived functor?
\end{question}

\section{Variations on a theme of Ado}

The Ado Theorem, one of the cornerstones of the modern theory of Lie algebras,
says that each finite-dimensional Lie algebra has a finite-dimensional 
faithful representation. Somewhat surprisingly, its proof is not that 
straightforward as one may expect for such a basic result: it involves universal
enveloping algebras -- infinite dimensional objects, and is strikingly different
for the cases of zero and positive characteristics. There exists a substantial 
body of literature with variants of the proof of the Ado theorem, but all of 
them follow, more or less, the same pattern. In \cite{ado} a different proof was
given, not appealing to the notion of universal enveloping algebra and intrinsic
to the category of finite-dimensional Lie algebras. Unfortunately, the proof is
valid for nilpotent Lie algebras and in characteristic zero only.

\begin{question}\label{quest-ado}
Give a characteristic-free, ``short'' and ``natural'' (i.e., not utilizing the 
notion of universal enveloping algebra or any other infinite-dimensional 
objects) proof of the full Ado Theorem.
\end{question}

In the standard proofs of the Ado Theorem, the case of nilpotent Lie algebras is
the most laborious part. Then, the general case is derived from the nilpotent 
one via the passage to the algebraic envelope, and employing a certain structure of
a faithful module built, again, with the help of universal enveloping algebra 
(and the PBW theorem).

To get a partial answer to Question \ref{quest-ado}, we may try to move along 
the same route, but employing ideas of \cite{ado}. As any finite-dimensional Lie
algebra is embedded into its algebraic envelope (first constructed by Malcev 
\cite{malcev}, and independently by Goto \cite{goto} and Matsushima 
\cite{matsushima}), it is enough to prove the Theorem for algebraic Lie 
algebras. In characteristic zero, the Levi--Malcev decomposition of any 
algebraic Lie algebra is of the form $L = S + T + N$ (direct sum of vector 
spaces), where $S$ is the semisimple part, $T$ is a torus acting on the 
nilradical $N$, and $[S,T] = 0$. As in the proof of the nilpotent case of the 
Theorem in \cite{ado}, it is enough to establish the existence of a faithful 
representation of $L$ not vanishing on any nonzero central element of $L$, and 
the latters lie in $N$. If, say, $N$ is $\mathbb N$-graded, then arguing like in \cite[Lemma 2.5]{ado}, we may 
construct a representation of $N$ in $L \otimes tK[t]/(t^n)$, for a suitable $n$, 
with required properties. Since $S$ and $T$ act on $N$ by derivations, we may 
extend this representation to the whole $L$. In this way we get a proof of the
Theorem for Lie
algebras whose algebraic envelope has an $\mathbb N$-graded nilpotent radical.

As for characteristic-free requirement, it is easy to see that all the 
reasonings of \cite{ado} remain valid over a field of characteristic $p$, if the
index of nilpotency of the Lie algebra in question is $< p$. But to give a full-blown characteristic-free proof will require, apparently, new ideas.

\section{Variations on a theme of Whitehead}

The Second Whitehead Lemma is another classical result saying that the second 
cohomology of a finite-dimensional semisimple Lie algebra of characteristic 
zero, with coefficients in arbitrary finite-dimensional module, vanishes.

In \cite{h2_is_zero} a curious ``almost converse'' was proved: a 
finite-dimensional Lie algebra of characteristic zero such that the second
cohomology in any its finite-dimensional module vanishes, is either semisimple,
or one-dimensional, or is the direct sum of a semisimple and one-dimensional 
algebra. According to \cite{makedonskii}, over an algebraically closed field
this is exactly the list of finite-dimensional Lie algebras having the tame 
representation type.

\begin{question}\label{q-c}
What is the reason that these two classes of Lie algebras coincide?
\end{question}

Note that in the modular case the situation is entirely different, due to the
result of Dzhumadil'daev \cite{dzhu-sbornik}: for any finite-dimensional Lie 
algebra over a field of positive characteristic, and any degree not exceeding 
the dimension of the algebra, there is a finite-dimensional module with 
non-vanishing cohomology in that degree. In particular, the one-dimensional Lie 
algebra is the only finite-dimensional Lie algebra with vanishing second 
cohomology in any finite-dimensional module.

\begin{question}
Is it true that over an (algebraically closed) field of positive characteristic,
the only finite-dimensional Lie algebra of tame representation type is 
one-dimensional?
\end{question}

If this question has an affirmative answer, then one may ask for a 
characteristic-free variant of Question \ref{q-c}; a satisfactory answer should
establish a bijection between these two classes of Lie algebras, without 
addressing peculiarities related to characteristic of the base field.

It is known that analogs of the Second Whitehead Lemma hold for other classes
of algebraic structures: Jordan algebras, alternative algebras, Lie triple 
systems, etc.

\begin{question}
Do some sort of converses hold for all these analogs of the Second Whitehead 
Lemma?
\end{question}

\section{Variations on a theme of Banach}

In \cite{banach} an attempt was made to trace the possible origins of a 
(vaguely formulated) question by Stefan Banach about ternary maps which are 
superpositions of binary maps.

This is possibly the most isolated topic in the present collection, though a 
nice relation with Lie theory exists: an answer to possible interpretations of 
the Banach question may be obtained using an idea borrowed from a pioneering 
paper by Jacobson about Lie triple systems.

In the process the following question arose. Let us count the number of ternary
maps $f: X \times X \times X \to X$ on a finite set $X$ of $n$ elements, which 
can be represented as a superposition of binary maps $*: X \times X \to X$. We 
get the following table:

\bigskip

\begin{center}
\begin{tabular}{|l|r|r|r|}
\hline
$n$ & $T_L(n)$ & $T_{LR}(n)$ & $T_{comm}(n)$ \\ \hline 
$1$ & $1$      & $1$         & $1$           \\ \hline 
$2$ & $14$     & $21$        & $5$           \\ \hline
$3$ & $19292$  & $38472$     & $48$          \\ \hline
\end{tabular}
\end{center}

\bigskip

Here $T_L(n)$ denotes the number of ternary maps represented in the form 
\begin{equation}\label{e-left}
f(x,y,z) = (x*y)*z
\end{equation}
for some binary map $*$, $T_{LR}(n)$ denotes the number of ternary maps 
represented either in the form (\ref{e-left}), or in the form 
\begin{equation}\label{e-right}
f(x,y,z) = x*(y*z) ,
\end{equation}
and $T_{comm}(n)$ denotes the number of ternary symmetric maps (i.e., invariant
under any permutation in $S_3$) represented in the form (\ref{e-left}). In the 
given range, the latter number coincides with the number of ternary symmetric 
maps represented in the form (\ref{e-left}) for some \emph{commutative} $*$.

As of time of this writing, the $3$-term sequences for $T_L(n)$ and $T_{LR}(n)$
were absent in The Online Encyclopedia of Integer Sequences, and among a dozen
or so sequences containing the $3$-term sequence for $T_{comm}(n)$, nothing 
seems to be relevant.

\begin{question}
Continue this table. Give formulas (closed form, or recurrent) for the numbers
$T_L(n)$, $T_{LR}(n)$, $T_{comm}(n)$.
\end{question}

\begin{question}
Is it true that for every $n$, any symmetric ternary map represented in the
form (\ref{e-left}) for some $*$, can  be represented in the form (\ref{e-left})
for some commutative $*$?
\end{question}

\section*{Acknowledgements}

The financial support of the Regional Authority of the 
Mo\-ra\-vi\-an-Si\-le\-si\-an Region (grant MSK 44/3316), and of the Ministry of
Education and Science of the Republic of Kazakhstan (grant 0828/GF4) is 
gratefully acknowledged.

\end{document}